\documentclass[12pt,a4paper]{article}
\usepackage{amsmath}
\usepackage{amsfonts}
\usepackage{amssymb}
\usepackage{amsmath,authblk,hyperref}
\newcommand{\RNum}[1]{\uppercase\expandafter{\romannumeral #1\relax}}

\title{\textbf{ Common Fixed Point Theorems for Self Mappings on a Complete Vector $S$-metric Space}}
\author[1]{Pooja Yadav\thanks{poojayadav.math.rs@igu.ac.in}}
\author[2]{Mamta Kamra\thanks{mkhaneja15@gmail.com}}
\author[3]{Rajpal\thanks{rajpal@rafflesuniversity.edu.in}}
\affil[1, 2]{Department of Mathematics, Indira Gandhi University  Meerpur(Rewari), Haryana-122502, India.}
\affil[3]{Department of Mathematics, Raffles University, Neemrana, Haryana, India.}
\date{}
\begin{document}
\maketitle
\begin{abstract}
\noindent$S$-metric space was introduced by Sedghi et al. \cite{Sedghi} in 2012. We derive some common fixed point results for self-mappings on vector valued complete $S$-metric space. In support of our results, we also give some examples. 
\end{abstract}
\section{Introduction}
Fixed point theory is amongst the crucial mathematical theory with applications in various branches of science. Banach contraction principle was derived first by S. Banach \cite{Banach} in 1922. It has a vital role in fixed point theory and became very famous due to iterations used in the theorem. The evaluation of fixed points of mappings satisfying many contractive conditions is at the center of research work and several vital results have been established by many authors. Over the last few years, several researchers have devoted themselves to define many variations of metric space. We give below some definitions and results which will help in proving our main results for vector $S$-metric spaces.\\ \\
\textbf{Definition 1.1}\cite{Kamra}
On a set $\complement$, a relation $\preceq$ is a partial order if it follows the conditions stated below:
\begin{itemize}
\item[(a)]$ \eta_1 \preceq\eta_1  $ \hspace{7.4cm} (reflexive)
\item[(b)]$\eta_1\preceq \eta_2  $ and $ \eta_2\preceq \eta_1  $ implies $\eta_1=\eta_2$                                                                                                                 \hspace{2.3cm} (anti-symmetry)
\item[(c)]$ \eta_1\preceq \eta_2  $ and $ \eta_2\preceq \eta_3  $\,\,implies $ \eta_1\preceq \eta_3$                                                                           
\hspace{2.4cm} (transitivity)
\end{itemize}
$\forall$ $\eta_1, \eta_2, \eta_3 \in \complement$  . The set $\complement$  with partial order $\preceq $ is known as partially ordered set (poset).\\
\vspace{.3cm}
A partially ordered set $(\complement, \preceq)$ is called  linearly ordered if for $ \eta_1, \eta_2  \in \complement$, we have either $\eta_1\preceq \eta_2$ or $\eta_2\preceq \eta_1$.\\ \\
\textbf{Definition 1.2}\cite{Kamra}
Let $\complement$ be  linear space which is real and $(\complement, \preceq)$ be a poset . Then the poset $(\complement, \preceq)$ is said to be an ordered linear space if it follows the  properties mentioned below:
\begin{itemize}
\item[(a)]$ p_1 \preceq p_2 \Longrightarrow p_1+p_3 \preceq p_2+p_3$
\item[(b)] $p_1 \preceq p_2 \Longrightarrow \omega p_1 \preceq \omega p_2$\hspace{1cm}$\forall p_1, p_2, p_3 \in \complement$ and $\omega >0$
\end{itemize} 
\textbf{Definition 1.3}\cite{Kamra}
A poset is called lattice if each set with two elements has an infimum and a supremum.\\ \\
\textbf{Definition 1.4}\cite{Kamra}
An ordered linear space where the ordering is lattice is called vector lattice. This is also called Riesz space.\\ \\
\textbf{Definition 1.5}\cite{Kamra}
A vector lattice $V$ is called Archimedean if\,\,  $inf\{\dfrac{1}{m}\vartheta\}=0$ for every  $\vartheta \in V^+$  where \[V^+= \{\vartheta \in V :\vartheta \geq 0\}.\] \\
\textbf{Definition 1.6}\cite{Cevik}
Let $V$ be a vector lattice and $\Re$ be a nonvoid set. A function $d:\Re \times\Re \rightarrow V$ is called vector metric on $\Re$ if it follows the conditions stated below:
\begin{itemize}
\item[(a)]$ d(\hslash_1, \hslash_2)=0$ iff $\hslash_1=\hslash_2$
\item[(b)]$ d(\hslash_1, \hslash_2) \preceq d(\hslash_1, \hslash_3) + d(\hslash_3, \hslash_2)   $\,\,\,\,$\forall \hslash_1, \hslash_2, \hslash_3\in \Re$
\end{itemize}
The triplet $(\Re, d, V)$ is called vector metric space.\\ \\
\textbf{Definition 1.7}\cite{Shahraki}
Let $\Re$ be a nonvoid set. A function $S:\Re\times \Re\times \Re\rightarrow [0, \infty)$ is called $S$-metric on $\Re$ if it follows the conditions stated below:
\begin{itemize}
\item [(a)] $S(\flat_1, \flat_2, \flat_3) \succeq 0$,
\item[(b)]$S(\flat_1, \flat_2, \flat_3) = 0$ iff $\flat_1 = \flat_2 = \flat_3$,
\item[(c)]$S(\flat_1, \flat_2, \flat_3) \preceq S(\flat_1, \flat_2, \alpha)+S(\flat_2, \flat_2, \alpha)+S(\flat_3, \flat_3, \alpha)$,
\end{itemize} 
\,\,for all  $\flat_1, \flat_2, \flat_3,\alpha \in \Re$.\\
The pair $(\Re, S)$ is called $S$-metric space .\\ \\
Now, we define vector valued $S$-metric space as follows:\\
\textbf{Definition 1.8}
Let $V$ be a vector lattice and $\Re$ be a nonvoid set. A function $S:\Re\times \Re\times \Re\rightarrow V$ is called vector $S$-metric on $\Re$ that satisfies the conditions mentioned below:
\begin{itemize}
\item [(a)] $S(\flat_1, \flat_2, \flat_3) \succeq 0$,
\item[(b)]$S(\flat_1, \flat_2, \flat_3) = 0$ iff $\flat_1 = \flat_2 = \flat_3$,
\item[(c)]$S(\flat_1, \flat_2, \flat_3) \preceq S(\flat_1, \flat_2, \alpha)+S(\flat_2, \flat_2, \alpha)+S(\flat_3, \flat_3, \alpha)$,
\end{itemize}
\,\,for all  $\flat_1, \flat_2, \flat_3,\alpha \in \Re$.\\
The triplet $(\Re, S, V)$ is called vector $S$-metric space.\\ \\
\textbf{Example 1.9}
Let $\Re$ be a nonvoid set and $V$ be a vector lattice. A function $S:V\times V\times V\rightarrow V$ is defined by 
\[S(\flat_1, \flat_2, \flat_3)=\vert(\flat_1, \flat_3)\vert + \vert(\flat_2, \flat_3)\vert \,\,\,\,\,\forall \flat_1, \flat_2, \flat_3 \in \Re \] then the triplet $(\Re, S, V)$ is vector $S$-metric space.\\\\
\textbf{Definition 1.10}
A sequence $\langle \vartheta_n \rangle$ in a vector $S$-metric space $(\Re, S, V)$ is called $V$-convergent to some $\vartheta \in V$ if there is a sequence $\langle \mu_n \rangle$ in $V$ satisfying $\mu_n \downarrow 0$ and $S(\vartheta_n, \vartheta_n, \vartheta) \leq \mu_n$ and denote it by $\vartheta_n \xrightarrow{S,V} \vartheta$.\\\\
\textbf{Definition 1.11}
A sequence $\langle \vartheta_n \rangle$ in a vector $S$-metric space $(\Re, S, V)$ is known as $V$-Cauchy sequence if $\exists \,\,\,\,\,\langle \mu_n \rangle \in V$ satisfying $\mu_n \downarrow 0$ and $S(\vartheta_n, \vartheta_n, \vartheta_{n+q}) \leq \mu_n$ \,\,\,$ \forall$ $q$ and $n$.\\\\
\textbf{Definition 1.12}
A vector $S$-metric space $(\Re, S, V)$ is called $V$-complete if all $V$-Cauchy sequence is $V$-convergent to a limit in $\Re$.\\ \\
\textbf{Lemma 1.13}\cite{Shahraki}
For a vector $S$-metric space $(\Re, S, V)$, \[S(\vartheta, \vartheta, \mu)=S(\mu, \mu, \vartheta) \,\,\,\,\,\,\,\forall \mu, \vartheta \in \Re.\] 
\textbf{Proof.} Using the condition (c) of definition (1.8), we have 
\begin{eqnarray}
S(\vartheta, \vartheta, \mu) &\leq &S(\vartheta, \vartheta, \vartheta)+S(\vartheta, \vartheta, \vartheta)+S(\mu, \mu, \vartheta)\\
&=&S(\mu, \mu, \vartheta) \nonumber\\
S(\mu, \mu, \vartheta)& \leq & S(\mu, \mu, \mu)+S(\mu, \mu, \mu)+S(\vartheta, \vartheta, \mu)\\
&=&S(\vartheta, \vartheta, \mu) \nonumber
\end{eqnarray}
By (1) and (2), we get $S(\vartheta, \vartheta, \mu)=S(\mu, \mu, \vartheta)$.\\
\section{Main Results}
\textbf{Lemma 2.1} Let   $(\Re, S, V)$ be a vector $S$-metric space which is complete and $V$-Archimedean. Let a sequence $\langle\hslash_\flat \rangle$ be in $\Re$ such that
\[ S(\hslash_\flat, \hslash_\flat, \hslash_{\flat+1}) \preceq {\alpha}S(\hslash_{\flat-1}, \hslash_{\flat-1}, \hslash_{\flat}) \,\,\,\,\,\, \forall \flat\in \mathbb{N} \hspace{1cm}(1)\]   
where $\alpha \in [0, 1)$. Then $\langle \hslash_\flat \rangle$  is a $V$- Cauchy sequence in $\Re$.\\
\textbf{Proof.} Using (1), we get
\[ S(\hslash_\flat, \hslash_\flat, \hslash_{\flat+1}) \preceq {\alpha}S(\hslash_{\flat-1}, \hslash_{\flat-1}, \hslash_{\flat}) \preceq \alpha^2S(\hslash_{\flat-2}, \hslash_{\flat-2}, \hslash_{\flat-1}) \preceq \dots \preceq \alpha^{\flat}S(\hslash_0, \hslash_0, \hslash_{1})  \] 
So, for ${\flat}>{\ell}$, we have 
\begin{eqnarray}
S(\hslash_\ell, \hslash_\ell, \hslash_{\flat}) & \preceq & 2S(\hslash_\ell, \hslash_\ell, \hslash_{\ell+1})+S(\hslash_{\flat}, \hslash_{\flat}, \hslash_{\ell+1}) \nonumber\\
& =& 2S(\hslash_\ell, \hslash_\ell, \hslash_{\ell+1})+S(\hslash_{\ell+1}, \hslash_{\ell+1}, \hslash_{\flat}) \nonumber\\
& \preceq & 2S(\hslash_\ell, \hslash_\ell, \hslash_{\ell+1})+2S(\hslash_{\ell+1}, \hslash_{\ell+1}, \hslash_{\ell+2})+ S(\hslash_{\flat}, \hslash_{\flat}, \hslash_{\ell+2})\nonumber\\
& = & 2S(\hslash_\ell, \hslash_\ell, \hslash_{\ell+1})+2S(\hslash_{\ell+1}, \hslash_{\ell+1}, \hslash_{\ell+2})+ S(\hslash_{\ell+2}, \hslash_{\ell+2}, \hslash_{\flat})\nonumber\\
& \preceq & 2S(\hslash_\ell, \hslash_\ell, \hslash_{\ell+1})+2S(\hslash_{\ell+1}, \hslash_{\ell+1}, \hslash_{\ell+2})+ \dots + S(\hslash_{\flat-1}, \hslash_{\flat-1}, \hslash_{\flat}) \nonumber\\
& \prec & 2S(\hslash_\ell, \hslash_\ell, \hslash_{\ell+1})+ 2S(\hslash_{\ell+1}, \hslash_{\ell+1}, \hslash_{\ell+2})+ \dots + 2S(\hslash_{\flat-1},\hslash_{\flat-1}, \hslash_{\flat}) \nonumber\\
& \prec & 2({\alpha}^{\ell}+{\alpha}^{\ell+1}+ \dots + {\alpha}^{\flat-1})S(\hslash_0, \hslash_0, \hslash_{1}) \nonumber\\
& \prec & 2{\alpha}^{\ell}(1+{\alpha}+{\alpha}^2+ \dots )S(\hslash_0, \hslash_0, \hslash_{1}) \nonumber\\
& \prec & 2\dfrac{{\alpha}^\ell}{1-{\alpha}}S(\hslash_0, \hslash_0, \hslash_{1})\downarrow 0\,\,\,\,\,\,\, \ell\rightarrow\infty. \nonumber
\end{eqnarray} 
Thus $\langle \hslash_\flat \rangle$ is a $V$-Cauchy sequence.\\ \\
\textbf{Theorem 2.2}
Let $(\Re, S, V)$ be a vector $S$-metric space which is complete and $V$- Archimedean. Let $K:\Re \rightarrow \Re$ be a continuous mapping and a map $f:\Re \rightarrow \Re$ which commutes with $K$. Suppose the conditions given below are satisfied;
\begin{itemize}
\item[(a)] $f(\Re) \subseteq K(\Re)$
\item[(b)]$S(f\hslash, f\hslash, f\vartheta) \preceq qU(\hslash, \hslash, \vartheta)$ for all $\hslash, \vartheta \in \Re$ where  $q \in \Big[0, \frac{1}{3}\Big)$ is a constant and \\
\begin{eqnarray}
U(\hslash, \hslash, \vartheta) &\in & \{S(K\hslash, K\hslash, K\vartheta), S(K\hslash, K\hslash, f\hslash), S(K\vartheta, K\vartheta, f\vartheta),\nonumber\\ & \,\,\,& S(K\hslash, K\hslash, f\vartheta), S(K\vartheta, K\vartheta, f\hslash)\}\nonumber
 \end{eqnarray}
\item[(c)]$K(\Re)$ or $ f(\Re) $  is $V$- complete as a subspace of $\Re$. Then, prove that $K$ and $f$  have a common fixed point which is unique .
\end{itemize} 
\textbf{Proof} Fix arbitrary $\vartheta_0 \in \Re$, so we can take sequence  $\langle \hslash_\flat \rangle$ in $\Re$ such that
\[\hslash_\flat =f\vartheta_\flat=K\vartheta_{\flat+1}\,\,\,\,\,\,\,\, {\flat}\geq 0.\] 
Then 
\begin{equation*}
S(\hslash_\flat, \hslash_\flat, \hslash_{\flat+1})=S(f\vartheta_\flat, f\vartheta_\flat, f\vartheta_{\flat+1}) \preceq qU(\vartheta_\flat, \vartheta_\flat, \vartheta_{\flat+1}) \hspace{2cm} (2) 
\end{equation*}
where 
\begin{eqnarray}
U(\vartheta_\flat, \vartheta_\flat, \vartheta_{\flat+1}) &\in & \{S(K\vartheta_{\flat}, K\vartheta_{\flat}, K\vartheta_{\flat+1}), S(K\vartheta_{\flat}, K\vartheta_{\flat}, f\vartheta_{\flat}), S(K\vartheta_{\flat+1},  \nonumber\\ &     \,\,\,& K\vartheta_{\flat+1}, f\vartheta_{\flat+1}), S(K\vartheta_{\flat}, K\vartheta_{\flat}, f\vartheta_{\flat+1}), S(K\vartheta_{\flat+1}, K\vartheta_{\flat+1}, f\vartheta_{\flat})\} \nonumber\\
&=& \{ S(\hslash_{\flat-1}, \hslash_{\flat-1}, \hslash_{\flat}), S(\hslash_{\flat-1}, \hslash_{\flat-1}, \hslash_{\flat}),S(\hslash_{\flat}, \hslash_{\flat}, \hslash_{\flat+1}),\nonumber\\ &     \,\,\,& S(\hslash_{\flat-1}, \hslash_{\flat-1}, \hslash_{\flat+1}), S(\hslash_{\flat}, \hslash_{\flat}, \hslash_{\flat})  \}\nonumber\\
&=& \{S(\hslash_{\flat-1}, \hslash_{\flat-1}, \hslash_{\flat}), S(\hslash_{\flat}, \hslash_{\flat}, \hslash_{\flat+1}),S(\hslash_{\flat-1}, \hslash_{\flat-1}, \hslash_{\flat+1}), 0\}\nonumber
\end{eqnarray}
The possible four cases are:\\
\begin{itemize}
\item[(i)]$S(\hslash_\flat, \hslash_\flat, \hslash_{\flat+1})  \preceq   qS(\hslash_{\flat-1}, \hslash_{\flat-1}, \hslash_{\flat})$
\item[(ii)]$ S(\hslash_\flat, \hslash_\flat, \hslash_{\flat+1})  \preceq   qS(\hslash_{\flat}, \hslash_{\flat}, \hslash_{\flat+1})$\\ 
and so\\ \[S(\hslash_{\flat}, \hslash_{\flat}, \hslash_{\flat+1})=0.\]
 
\item[(iii)]$S(\hslash_\flat, \hslash_\flat, \hslash_{\flat+1})  \preceq qS(\hslash_{\flat-1}, \hslash_{\flat-1}, \hslash_{\flat+1}) \preceq 2qS(\hslash_{\flat-1}, \hslash_{\flat-1}, \hslash_{\flat})+qS(\hslash_{\flat+1}, \hslash_{\flat+1}, \hslash_{\flat}).$\\ and \[S(\hslash_\flat, \hslash_\flat, \hslash_{\flat+1}) \preceq \dfrac{2q}{(1-q)} S(\hslash_{\flat-1}, \hslash_{\flat-1}, \hslash_{\flat})\]
\item[(iv)]$S(\hslash_\flat, \hslash_\flat, \hslash_{\flat+1})  \preceq   q.0=0$\\
and so\\ \[S(\hslash_\flat, \hslash_\flat, \hslash_{\flat+1})=0.\]\\
\end{itemize}
Thus $S(\hslash_\flat, \hslash_\flat, \hslash_{\flat+1})  \preceq   \sigma S(\hslash_{\flat-1}, \hslash_{\flat-1}, \hslash_{\flat})$\,\,\,\, where $\sigma \in \Big\{q, \dfrac{2q}{(1-q)}\Big\}<1$ \\ \\
Since $V$ is Archimedean, by lemma (2.1) $\langle \hslash_\flat \rangle$ is a $V$-Cauchy sequence and  range of $f$ is contained in the range of $K$ and atleast one range is $V$-complete, there exist  $\hslash \in K(\Re)$ such that $K\vartheta_{\flat+1} \xrightarrow {S, V} \hslash$. Hence there exist  a sequence $\langle \alpha_\flat \rangle$ in $V$ 
such that $\alpha_\flat \downarrow 0$ and 
\[S(K\vartheta_\flat, K\vartheta_\flat, \hslash) \preceq \alpha_\flat.\]
Thus \begin{eqnarray}
\hslash_\flat=f\vartheta_\flat=K\vartheta_{\flat+1} \xrightarrow {S, V} \hslash.
\end{eqnarray}
We prove that 
\[K\hslash= f\hslash= \hslash \]
Now,
\begin{eqnarray}
S(K\hslash, K\hslash, f\hslash) &\preceq& 2S(K\hslash, K\hslash, fK\vartheta_\flat)+S(f\hslash, f\hslash, fK\vartheta_\flat)\nonumber\\
 S(K\hslash, K\hslash, f\hslash)  &\preceq& 2S(K\hslash, K\hslash, fK\vartheta_\flat)+S(fK\vartheta_\flat, fK\vartheta_\flat, f\hslash)
\end{eqnarray}
Also, we have
\begin{eqnarray} 
S(fK\vartheta_\flat, fK\vartheta_\flat, f\hslash)\, \preceq \,qU(K\vartheta_\flat, K\vartheta_\flat, \hslash) \nonumber
\end{eqnarray}
Then (3) becomes\\
\begin{eqnarray}
 S(K\hslash, K\hslash, f\hslash)  &\preceq& 2S(K\hslash, K\hslash, fK\vartheta_\flat)+qU(K\vartheta_\flat, K\vartheta_\flat, \hslash)
\end{eqnarray}
where \\
\begin{eqnarray}
U(K\vartheta_\flat, K\vartheta_\flat, \hslash)  &\in &  \{ S(K^2\vartheta_\flat, K^2\vartheta_\flat, K\hslash), S(K^2\vartheta_\flat, K^2\vartheta_\flat, fK\vartheta_\flat),\nonumber\\ &     \,\,\,& \hspace{.1cm}S(K\hslash, K\hslash, f\hslash),   S(K^2\vartheta_\flat, K^2\vartheta_\flat, f\hslash),\nonumber\\ &     \,\,\,\,\,\,\,\,& S(K\hslash, K\hslash, fK\vartheta_\flat)\}
\end{eqnarray}
 Since $f$ commutes with $K$ and by using continuity of $K$, we get
\[fK\vartheta_\flat=Kf\vartheta_\flat \xrightarrow{S, V} K\hslash \] and by using (3)\[K^2\vartheta_{\flat+1} \xrightarrow{S, V} K\hslash,\]
then there exist a sequence $\langle \alpha_\flat \rangle$ and $\langle \beta_\flat \rangle$ in $V$ such that $\alpha_\flat \downarrow 0$ and $\beta_\flat \downarrow 0$ , then we have 
\begin{eqnarray}
S(K\hslash, K\hslash, fK\vartheta_\flat) &\preceq & \alpha_\flat \nonumber
\end{eqnarray}
and
\begin{eqnarray}
 S(K^2\vartheta_\flat, K^2\vartheta_\flat, K\hslash) &\preceq & \beta_\flat. \nonumber
\end{eqnarray}
By (5) and (6), we have the following cases:\\
\begin{eqnarray}
(i)\hspace{2cm} S(K\hslash, K\hslash, f\hslash) &\preceq & 2S(K\hslash, K\hslash, fK\vartheta_\flat) + qS(K^2\vartheta_\flat, K^2\vartheta_\flat, K\hslash)\nonumber\\
&\preceq& 2\alpha_\flat +q\beta_\flat \nonumber\\
(ii)\hspace{2cm}S(K\hslash, K\hslash, f\hslash) &\preceq & 2S(K\hslash, K\hslash, fK\vartheta_\flat) + qS(K^2\vartheta_\flat, K^2\vartheta_\flat, fK\vartheta_\flat) \nonumber\\
&\preceq &  2S(K\hslash, K\hslash, fK\vartheta_\flat)  +q[2S(K^2\vartheta_\flat, K^2\vartheta_\flat, K\hslash)+\nonumber\\ &     \,\,\,& S(fK\vartheta_\flat, fK\vartheta_\flat, K\hslash)]\nonumber\\
&= & 2S(K\hslash, K\hslash, fK\vartheta_\flat) +2qS(K^2\vartheta_\flat, K^2\vartheta_\flat, K\hslash)+\nonumber\\ &     \,\,\,& qS(K\hslash, K\hslash, fK\vartheta_\flat) \nonumber\\
&\preceq& (2+q)\alpha_\flat + 2q\beta_\flat\nonumber\\
(iii)\hspace{2cm}  S(K\hslash, K\hslash, f\hslash) &\preceq & 2S(K\hslash, K\hslash, fK\vartheta_\flat)  + q S(K\hslash, K\hslash, f\hslash)\nonumber\\
\,\,\,\,\,\,(1-q)S(K\hslash, K\hslash, f\hslash) &\preceq &  2\alpha_\flat\nonumber\\
S(K\hslash, K\hslash, f\hslash)&\preceq& \dfrac{2\alpha_\flat }{(1-q)} \nonumber
\end{eqnarray}
\begin{eqnarray}
(iv)\hspace{2cm}S(K\hslash, K\hslash, f\hslash) &\preceq & 2S(K\hslash, K\hslash, fK\vartheta_\flat) + qS(K^2\vartheta_\flat, K^2\vartheta_\flat, f\hslash)\nonumber\\
&\preceq & 2S(K\hslash, K\hslash, fK\vartheta_\flat) +2qS(K^2\vartheta_\flat, K^2\vartheta_\flat, K\hslash)+\nonumber\\ &     \,\,\,& qS(f\hslash, f\hslash, K\hslash)\nonumber\\
&= & 2S(K\hslash, K\hslash, fK\vartheta_\flat) +2qS(K^2\vartheta_\flat, K^2\vartheta_\flat, K\hslash)+\nonumber\\ &     \,\,\,& qS(K\hslash, K\hslash, f\hslash)\nonumber\\
(1-q)S(K\hslash, K\hslash, f\hslash) &\preceq &  2\alpha_\flat +2q\beta_\flat\nonumber\\
S(K\hslash, K\hslash, f\hslash) &\preceq & \dfrac{2\alpha_\flat +2q\beta_\flat}{(1-q)}\nonumber\\
(v)\hspace{2cm}S(K\hslash, K\hslash, f\hslash) &\preceq & 2S(K\hslash, K\hslash, fK\vartheta_\flat) + qS(K\hslash, K\hslash, fK\vartheta_\flat) \nonumber\\
&\preceq& (2+q)S(K\hslash, K\hslash, fK\vartheta_\flat)\nonumber\\
&\preceq& (2+q)\alpha_\flat \nonumber\\
&\preceq & 3\alpha_\flat \nonumber
\end{eqnarray}\\
In the last inequality of each case, the infimum on the right hand side is 0. So we get\\
 \[S(K\hslash, K\hslash, f\hslash)=0.\]
This implies \,\,\,\,\,\begin{eqnarray}
K\hslash=f\hslash
\end{eqnarray}
So\\
\begin{eqnarray}
S(K\hslash, K\hslash, \hslash) & \preceq& 2S(K\hslash, K\hslash, f\vartheta_\flat)+ S(\hslash, \hslash, f\vartheta_\flat)  \nonumber\\
 & =& 2S(K\hslash, K\hslash, f\vartheta_\flat)+ S(f\vartheta_\flat, f\vartheta_\flat, \hslash)  \nonumber\\
 & \preceq&  S(f\vartheta_\flat, f\vartheta_\flat, \hslash)+ 2S(f\hslash, f\hslash, f\vartheta_\flat) \nonumber\\
  & =&  S(f\vartheta_\flat, f\vartheta_\flat, \hslash)+ 2S(f\vartheta_\flat, f\vartheta_\flat, f\hslash) \nonumber\\
S(K\hslash, K\hslash, \hslash) &\preceq & S(f\vartheta_\flat, f\vartheta_\flat, \hslash) + 2qU(\vartheta_\flat, \vartheta_\flat, \hslash)
\end{eqnarray} 
where \\
\begin{eqnarray}
U(\vartheta_\flat, \vartheta_\flat, \hslash) & \in & \{S(K\vartheta_\flat, K\vartheta_\flat, K\hslash), S(K\vartheta_\flat, K\vartheta_\flat, f\vartheta_\flat), S(K\hslash, K\hslash, f\hslash), \nonumber\\ &     \,\,\,& S(K\vartheta_\flat, K\vartheta_\flat, f\hslash), S(K\hslash, K\hslash, f\vartheta_\flat)\}\nonumber\\
&=& \{S(K\vartheta_\flat, K\vartheta_\flat, K\hslash), S(K\vartheta_\flat, K\vartheta_\flat, f\vartheta_\flat), 0, S(K\hslash, K\hslash, f\vartheta_\flat)\}
\nonumber\\ &     \,\,\,& \hspace{9cm}
\end{eqnarray}
there exist a sequence $\langle d_n\rangle$ and $\langle g_n\rangle$ in $V$ such that $d_n \downarrow 0$ and $g_n \downarrow 0$, then we have \\
\begin{eqnarray}
S(K\vartheta_\flat, K\vartheta_\flat, \hslash) &\preceq & d_n  \nonumber
\end{eqnarray}
and
\begin{eqnarray}
S(f\vartheta_\flat, f\vartheta_\flat, \hslash) &\preceq & g_n .\nonumber
\end{eqnarray}
Using (8) and (9), we have the following cases:
\begin{eqnarray}
(i) \hspace{2cm} S(K\hslash, K\hslash, \hslash) & \preceq &  S(f\vartheta_\flat, f\vartheta_\flat, \hslash) + 2qS(K\vartheta_\flat, K\vartheta_\flat, K\hslash) \nonumber\\
& \preceq & S(f\vartheta_\flat, f\vartheta_\flat, \hslash)  + 2q[2S(K\vartheta_\flat, K\vartheta_\flat, \hslash)+\nonumber\\ &\,\,& S(K\hslash, K\hslash, \hslash) ] \nonumber\\
(1-2q) S(K\hslash, K\hslash, \hslash) & \preceq & g_n+ 4qd_n \nonumber\\
 S(K\hslash, K\hslash, \hslash) & \preceq& \dfrac{g_n+ 4qd_n }{(1-2q)} \nonumber
\\
(ii) \hspace{2cm}  S(K\hslash, K\hslash, \hslash) & \preceq &  S(f\vartheta_\flat, f\vartheta_\flat, \hslash) + 2qS(K\vartheta_\flat, K\vartheta_\flat, f\vartheta_\flat) \nonumber\\
& \preceq &  S(f\vartheta_\flat, f\vartheta_\flat, \hslash) + 2q[2S(K\vartheta_\flat, K\vartheta_\flat, \hslash)+\nonumber\\ &\,\,& S(f\vartheta_\flat, f\vartheta_\flat, \hslash) ] \nonumber\\
& = &  S(f\vartheta_\flat, f\vartheta_\flat, \hslash) + 2q[2S(K\vartheta_\flat, K\vartheta_\flat, \hslash)+ \nonumber\\ &\,\,&S(\hslash, \hslash, f\vartheta_\flat) ] \nonumber\\
S(K\hslash, K\hslash, \hslash) & \preceq &  g_n + 2q(2d_n + g_n) \nonumber\\
& \preceq & g_n (1+2q)+ 4qd_n \nonumber
 \end{eqnarray}
\begin{eqnarray} 
(iii) \hspace{2cm} S(K\hslash, K\hslash, \hslash) & \preceq &   S(f\vartheta_\flat, f\vartheta_\flat, \hslash)+ 2q.0 \nonumber\\
S(K\hslash, K\hslash, \hslash) & \preceq &  g_n \nonumber\\
(iv) \hspace{2cm} S(K\hslash, K\hslash, \hslash) & \preceq &  S(f\vartheta_\flat, f\vartheta_\flat, \hslash)+ 2qS(K\hslash, K\hslash, f\vartheta_\flat)  \nonumber\\
& \preceq&  S(f\vartheta_\flat, f\vartheta_\flat, \hslash) + 2q[2S(K\hslash, K\hslash, \hslash)+ \nonumber\\ &\,\,&S(f\vartheta_\flat, f\vartheta_\flat, \hslash)] \nonumber\\
& =&  S(f\vartheta_\flat, f\vartheta_\flat, \hslash) + 2q[2S(K\hslash, K\hslash, \hslash)+ \nonumber\\ &\,\,&S(\hslash, \hslash, f\vartheta_\flat)] \nonumber\\
(1-4q)S(K\hslash, K\hslash, \hslash) & \preceq & g_n+ 2qg_n \nonumber\\
S(K\hslash, K\hslash, \hslash) & \preceq &   \dfrac{(1+2q)}{(1-4q)} g_n. \nonumber
\end{eqnarray}
In the last inequality of each case, the infimum on the right hand side is 0. So we get\\
\[S(K\hslash, K\hslash, \hslash)=0  \]
So \[K\hslash=\hslash.\]
From (7) \[ K\hslash= f\hslash= \hslash\]
So, $f$ and $K$ have common fixed  point $\hslash$  .\\
If  $K$ and $f$ have another common fixed  point $\mu_1$ then 
 \[ K\mu_1= f\mu_1= \mu_1.\]
From hypothesis (b)\[S(\hslash, \hslash, \mu_1)=S(f\hslash, f\hslash, f\mu_1) \preceq qU(\hslash, \hslash, \mu_1)\]
 where\\
\begin{eqnarray}
U(\hslash, \hslash, \mu_1) & \in & \{S(K\hslash, K\hslash, K\mu_1), S(K\hslash, K\hslash, f\hslash), S(K\mu_1, K\mu_1, f\mu_1 ),\nonumber\\ & \,\,\,& S(K\hslash, K\hslash, f\mu_1), S(K\mu_1, K\mu_1, f\hslash)\} \nonumber\\
&=& \{0, S(\hslash, \hslash, \mu_1).\} \nonumber
\end{eqnarray}
Thus \[ S(\hslash, \hslash, \mu_1)=0.\]
So \[\hslash=\mu_1.\]
Hence $f$ and  $K$ have a common fixed point $\hslash$ that is unique. \\ \\
\textbf{Corollary 2.3}
Let $(\Re, S, V)$ be a vector $S$-metric space which is complete and $V$- Archimedean. Let $K:\Re \rightarrow \Re$ be a map which is continuous and $f:\Re \rightarrow \Re$ be a map  which commutes with $K$. Also let $f$ and $K$ satisfy $f(\Re) \subseteq K(\Re)$ and   $S(f\hslash, f\hslash, f\vartheta)\preceq qS(K\hslash, K\hslash, K\vartheta)$ for all $q\in \Big[0,\dfrac{1}{3}\Big)$ and  $\hslash, \vartheta\in \Re$.Then $K$ and $f$  have  common fixed point which is unique.\\ \\
\textbf{Theorem 2.4}
Let $(\Re, S, V)$ be a vector $S$-metric space which is complete and $V$-Archimedean. Suppose a map $K^2:\Re \rightarrow \Re$ is continuous and a map  $f:\Re \rightarrow \Re$  which commutes with $K$. Suppose the conditions given below are satisfied:
\begin{itemize}
\item[(i)]$fK(\Re) \subseteq K^2(\Re) $
\item[(ii)]$S(f\hslash, f\hslash, f\vartheta) \preceq qU(\hslash, \hslash, \vartheta)$ for all $\hslash, \vartheta \in \Re$ where $q \in \Big[0, \dfrac{1}{3}\Big)$ is a constant and\\
 \begin{eqnarray} 
U(\hslash, \hslash, \vartheta) & \in & \{S(K\hslash, K\hslash, K\vartheta), S(K\hslash, K\hslash, f\hslash), S(K\vartheta, K\vartheta, f\vartheta),  \nonumber\\ & \,\,\,& \dfrac{1}{3}[S(K\hslash, K\hslash, f\vartheta)+ S(K\vartheta, K\vartheta, f\hslash)]\}\nonumber
\end{eqnarray}
\item[(iii)]$K(\Re)$ or $f(\Re)$ as a subspace of $\Re$ which is $V$-complete.\\
Then prove that $K$ and $f$ have a common fixed point which is unique.
\end{itemize}
\textbf{Proof.} Fix arbitrary $\vartheta_0 \in K(\Re)$, so we can take sequence  $\langle \hslash_\flat \rangle$ in $K(\Re)$ such that
\[\hslash_\flat =f\vartheta_\flat=K\vartheta_{\flat+1}\,\,\,\,\,\,\,\, \flat \geq 0.\] 
Now \[K\hslash_{\flat+1}=Kf\vartheta_{\flat+1}=fK\vartheta_{\flat+1}=f\hslash_\flat=p_\flat \,\,\,\,\,\,\,\,\, \flat\geq 0\]
It can be shown as in theorem(2.2) that $\langle p_\flat \rangle$ is a $V$-Cauchy sequence and converge to  $p \in \Re$ and 
\[K^2 p=fKp\]\\
\begin{eqnarray}
\lim_{\flat\rightarrow\infty}Kf\vartheta_\flat= \lim_{\flat\rightarrow\infty}K\hslash_{\flat}= \lim_{\flat\rightarrow\infty}p_{\flat-1}=p. 
\end{eqnarray}
It follows that \\
\begin{eqnarray}
\lim_{\flat\rightarrow\infty}K^4 \vartheta_\flat& = &\lim_{\flat\rightarrow\infty}K^3(K\vartheta_\flat) \nonumber\\
& = & \lim_{\flat\rightarrow\infty}K^3(\hslash_{\flat-1}) \nonumber\\
& = & \lim_{\flat\rightarrow\infty}K^2(K\hslash_{\flat-1}) \nonumber\\
& = & \lim_{\flat\rightarrow\infty}K^2(p_{\flat-2}) \nonumber\\
& = & K^2p, \,\,\,because\,\, K^2 \,\,is\,\, continuous. 
\end{eqnarray}
 So we get\\
\begin{eqnarray}
S(K^2 p, K^2 p, fKp) &\preceq & 2S(K^2 p, K^2 p, fK^3 \vartheta_\flat)+S(fKp, fKp, fK^3 \vartheta_\flat ) \nonumber\\ 
&= & 2S(K^2 p, K^2 p, fK^3 \vartheta_\flat)+S(fK^3 \vartheta_\flat, fK^3 \vartheta_\flat, fKp ) \nonumber\\
S(K^2 p, K^2 p, fKp)  &\preceq & 2S(K^2 p, K^2 p, fK^3 \vartheta_\flat)+ qU(K^3 \vartheta_\flat, K^3 \vartheta_\flat, Kp)  
\end{eqnarray}
where\\
\begin{eqnarray}
U(K^3 \vartheta_\flat, K^3 \vartheta_\flat, Kp) &\in & \{S(K^4 \vartheta_\flat, K^4 \vartheta_\flat, K^2p), S(K^4 \vartheta_\flat, K^4 \vartheta_\flat, fK^3 \vartheta_\flat),\nonumber\\ & \,\,\,& S(K^2 p, K^2 p, fKp), \dfrac{1}{3}[S(K^4 \vartheta_\flat, K^4 \vartheta_\flat, fKp)+\nonumber\\ & \,\,\,&S(K^2 p, K^2 p, fK^3 \vartheta_\flat)]\} 
\end{eqnarray}
Since $K^2$ is continuous and by using (10), we get
\[K^3 f\vartheta_\flat \xrightarrow{S, V} K^2p\] and \[K^4 \vartheta_\flat \xrightarrow{S, V} K^2p\]
then there exist a sequence $\langle \alpha_\flat \rangle$ and $\langle \beta_\flat \rangle$ in $V$ such that $\alpha_\flat \downarrow 0$ and $\beta_\flat \downarrow 0$, then we have 
\begin{eqnarray}
\hspace{1.5cm}S(K^2p, K^2 p, fK^3 \vartheta_\flat)&\preceq& \alpha_\flat \nonumber\\ S(K^4 \vartheta_\flat, K^4 \vartheta_\flat, K^2p)&\preceq& \beta_\flat.\nonumber
\end{eqnarray}
Using (12) and (13), then we have the following:
\begin{eqnarray}
Case(i) \hspace{1.5cm} S(K^2 p, K^2 p, fKp) &\preceq & 2S(K^2 p, K^2 p, fK^3 \vartheta_\flat)+q S(K^4 \vartheta_\flat, K^4 \vartheta_\flat, K^2p) \nonumber\\
&\preceq & 2\alpha_\flat  + q \beta_\flat\nonumber \\
Case(ii) \hspace{1.5cm} S(K^2 p, K^2 p, fKp) &\preceq & 2S(K^2 p, K^2 p, fK^3 \vartheta_\flat)+q   S(K^4 \vartheta_\flat, K^4 \vartheta_\flat, fK^3 \vartheta_\flat) \nonumber\\
 &\preceq & 2S(K^2 p, K^2 p, fK^3 \vartheta_\flat)+q[2S(K^4 \vartheta_\flat, K^4 \vartheta_\flat, K^2p)+\nonumber\\ & \,\,\,&S(fK^3 \vartheta_\flat, fK^3 \vartheta_\flat,  K^2 p)]\nonumber\\
 &= & 2S(K^2 p, K^2 p, fK^3 \vartheta_\flat)+q[2S(K^4 \vartheta_\flat, K^4 \vartheta_\flat, K^2p)+\nonumber\\ & \,\,\,&S(K^2p, K^2p,  K^3 f\vartheta_\flat)]\nonumber\\
&\preceq & 2\alpha_\flat+q(2\beta_\flat+\alpha_\flat)\nonumber\\
&\preceq &(2+q)\alpha_\flat+2q\beta_\flat\nonumber\\
Case(iii) \hspace{1.5cm}S(K^2 p, K^2 p, fKp) &\preceq & 2S(K^2 p, K^2 p, fK^3 \vartheta_\flat)+q S(K^2 p, K^2 p, fKp)\nonumber\\
(1-q)S(K^2 p, K^2 p, fKp) &\preceq & 2\alpha_\flat \nonumber\\
S(K^2 p, K^2 p, fKp) &\preceq &\dfrac{ 2\alpha_\flat }{(1-q)}\nonumber\\
Case(iv) \hspace{1.5cm} S(K^2 p, K^2 p, fKp) &\preceq & 2S(K^2 p, K^2 p, fK^3 \vartheta_\flat)+\dfrac{q}{3}[S(K^4 \vartheta_\flat, K^4 \vartheta_\flat, fKp)+\nonumber\\ & \,\,\,&S(K^2 p, K^2 p, fK^3 \vartheta_\flat)]\nonumber\\
S(K^2 p, K^2 p, fKp) &\preceq & 2S(K^2 p, K^2 p, fK^3 \vartheta_\flat)+\dfrac{q}{3}[2S(K^4 \vartheta_\flat, K^4 \vartheta_\flat, K^2 p)+\nonumber\\ & \,\,\,&S(fKp, fKp, K^2 p)+S(K^2 p, K^2 p, fK^3 \vartheta_\flat)]\nonumber
\end{eqnarray}
\begin{eqnarray} 
&= & 2S(K^2 p, K^2 p, fK^3 \vartheta_\flat)+\dfrac{q}{3}[2S(K^4 \vartheta_\flat, K^4 \vartheta_\flat, K^2 p)+\nonumber\\ & \,\,\,&S(K^2 p, K^2 p, fKp)+S(K^2 p, K^2 p, fK^3 \vartheta_\flat)]\nonumber\\
(1-\dfrac{q}{3})S(K^2 p, K^2 p, fKp)&\preceq & 2\alpha_\flat + \dfrac{q}{3}[2\beta_\flat  + \alpha_\flat ]\nonumber\\
(3-q)S(K^2 p, K^2 p, fKp)&\preceq & 6\alpha_\flat  + q[2\beta_\flat  + \alpha_\flat ]\nonumber\\
S(K^2 p, K^2 p, fKp)&\preceq & \dfrac{(6+q)\alpha_\flat +2q\beta_\flat }{(3-q)}\nonumber
\end{eqnarray}
In the last inequality of each case, the infimum on the right hand side is 0. So \[S(K^2 p, K^2 p, fKp)=0\]
that is\[K^2 p= fKp\]
Putting in $S(f\hslash, f\hslash, f\vartheta) \preceq qU(\hslash, \hslash, \vartheta)$, $\hslash=fKp, \vartheta=Kp$ then we get $f(fKp)=fK$, we get $f(fKp)= fKp$. Because $K^2 p= fKp$ i.e. $K(Kp)=f(Kp)$, we have $K(fKp)=fK^2 p=f(fKp)= fKp$. So  $f$ and $K$ have a common fixed point $fKp$ .\\
If  $f$ and $K$ have a common fixed point $f_1K_1p_1$ then 
\[K(f_1K_1p_1)=f(f_1K_1p_1)=f_1K_1p_1\]
Now \[S(fKp, fKp, f_1K_1p_1)=S(f(fKp), f(fKp), f(f_1K_1p_1))\preceq qU(fKp, fKp, f_1K_1p_1 )\]
where\\
\begin{eqnarray}
U(fKp, fKp, f_1K_1p_1 ) &\in& \{S(K(fKp), K(fKp), K(f_1K_1p_1)), S(K(fKp), K(fKp), f(fKp)), \nonumber\\ & \,\,\,& S(K(f_1K_1p_1), K(f_1K_1p_1), f(f_1K_1p_1)), S(K(fKp), K(fKp), f(f_1K_1p_1))), \nonumber\\ & \,\,\,&S(K(f_1K_1p_1), K(f_1K_1p_1), f(fKp))\}\nonumber\\
&\in&\{0, S(fKp, fKp, f_1K_1p_1)\}
\end{eqnarray}
Thus\\
\[S(fKp, fKp, f_1K_1p_1)=0.\]
So
\[fKp=f_1K_1p_1.\]
Hence $f$ and $K$ have a common fixed point $fKp$ which is unique.\\ \\
\textbf{Example 2.5}
Let $V=\mathbb{R}^2_+$ with coordinatewise ordering and $\Re=\mathbb{R}$
\[S(\hslash, \vartheta, \eta)=(\rho\vert \hslash-\eta \vert, \sigma\vert \vartheta-\eta \vert)\]
where $\rho, \sigma > 0$, and $\hslash, \vartheta, \eta \in \Re$. Then 
\[f\hslash=\hslash^2+5\]and\[K\hslash=2\hslash^2.\]
We have \[S(f\hslash, f\vartheta, f\eta )=(\rho\vert \hslash^2-\eta^2 \vert,\sigma\vert \vartheta^2-\eta^2 \vert)=\dfrac{1}{2}S(K\hslash, K\vartheta, K\eta )\leq qS(K\hslash, K\vartheta, K\eta ) \]
for $q\in [\frac{1}{2}, 1)$, $f(\Re)=[5,\infty)\subseteq [0, \infty)=K(\Re)$ and  a self map $K $ is  continuous on $\Re$ and $f(\Re)$ is $V$-complete subspace of $\Re$. Hence $K$ and $f$ have common fixed point that is unique.\\ \\ 
\textbf{Example 2.6}
Let $V=\mathbb{R}$ with coordinatewise ordering and $\Re=[0, 1]$
\[S(\hslash, \vartheta, \eta)=\vert \hslash-\eta \vert+\vert \vartheta-\eta \vert\]
where $\hslash, \vartheta, \eta \in \Re$. Then 
\[f\hslash=\dfrac{\hslash}{4}\]and\[K\hslash=\dfrac{\hslash}{2} \]
We have \[S(f\hslash, f\vartheta, f\eta )=\Big\vert \dfrac{\hslash}{4}-\dfrac{\eta}{4} \Big\vert+\Big\vert \dfrac{\vartheta}{4}-\dfrac{\eta}{4} \Big\vert=\dfrac{1}{2}S(K\hslash, K\vartheta, K\eta )\leq qS(K\hslash, K\vartheta, K\eta ) \]
for $q\in [\frac{1}{2}, 1)$, $f(\Re)=\Big[0, \dfrac{1}{4}\Big]\subseteq \Big[0, \dfrac{1}{2}\Big]=K(\Re)$ and  a self map $K $ is  continuous on $\Re$ and $f(\Re)$ is $V$-complete subspace of $\Re$. Hence $K$ and $f$ have common fixed point that is unique.\\ \\


\begin{thebibliography}{}
\bibitem{Altun} Altun I. , Cevik C. , Some common fixed point theorems in vector metric spaces, Filomat, 25(1) (2011) 105-113.
\bibitem{Banach} Banach S., Sur les operations dans les ensembles abstraits el leur application aux equations integrals, Fund. Math., 3 (1992) 133–181.

\bibitem{Aliprantis} C. D. Aliprantis, K. C. Border, \emph{Infinite Dimensional Analysis}, Verlag, Berlin, 1999.
\bibitem{Cevik} Cevik C., Altun I., \emph{Vector metric spaces and some properties}, Topal. Met. Nonlin. Anal., 34(2) (2009) 375-382.

\bibitem{Kamra} Kamra M., Kumar S., Sarita K., \emph{Some fixed point theorems for self mappings on vector b-metric spaces}, Global Journal of Pure and Applied Mathematics, 14(11) (2018) 1489-1507 .
\bibitem{Sedghi} Sedghi S., Shobe N. , Aliouche A., \emph{A generalization of fixed point theorem in S-metric spaces}, Mat. Vesnik, 64 (2012) 258-266.
\bibitem{Shahraki} Shahraki M., Sedghi S., Aleomraninejad S. M. A., Mitrovic Z. D., \emph{Some fixed point results on S-metric spaces}, Acta Univ. Sapientiae, Mathematica, 12(2) (2020) 347-357.



\end{thebibliography}
\end{document}